\documentclass[a4paper, 11pt]{amsart}
\usepackage{amsmath,amsthm,amsfonts,amssymb}
\theoremstyle{plain}

\newtheorem{teo}{Theorem}[section]
\newtheorem{lem}[teo]{Lemma}

\theoremstyle{definition}
\newtheorem*{conA}{Conjecture A}
\newtheorem*{thmB}{Theorem B}
\newtheorem*{thmC}{Theorem C}
\newtheorem*{thmD}{Theorem D}
\newtheorem{step}{Step}

\newcommand{\Ker}{\operatorname{Ker}}
\newcommand{\Lin}{\operatorname{Lin}}

\newcommand{\Irr}{\operatorname{Irr}} 
 
\newcommand{\Z}{\mathbf{Z}}
\newcommand{\C}{\mathbf{C}}

\newcommand{\nor}{\vartriangleleft}
\newcommand{\sbs}{\subseteq}

\title{Homogeneous products of characters}

\author{}

\date{}

\begin{document}

\maketitle

\centerline{by}

\centerline{Edith Adan-Bante}
\centerline{Department of Mathematics, University of Illinois at
Urbana-Champaign}
\centerline{Urbana IL 61801 USA}
\centerline{E-mail: adanbant@uiuc.edu}
\centerline{Maria Loukaki}
\centerline{Department of Applied Mathematics, University of Crete}
\centerline{Knosou Av. 74 GR-71409}
\centerline{Heraklion-Crete GREECE}
\centerline{E-mail: loukaki@tem.uoc.gr}
\centerline{and}
\centerline{Alexander Moret\'o}
\centerline{Departament d'\`Algebra, Universitat de Val\`encia}
\centerline{46100 Burjassot. Val\`encia SPAIN}
\centerline{E-mail: mtbmoqua@lg.ehu.es}

\noindent

The first author was partially supported by the 
National Science Foundation by
grant DMS 99-70030.
The third author was supported by the Basque Government and the Spanish
Ministerio de Ciencia y Tecnolog\'{\i}a.

\newpage

\section{Introduction}

I. M. Isaacs \cite{isa} has conjectured that if the product of two faithful
irreducible characters of a solvable group is irreducible, then the group is
cyclic. 
In this note we discuss the following conjecture, which generalizes Isaacs
conjecture.

\begin{conA}
Suppose that $G$ is solvable and that $\psi,\varphi\in\Irr(G)$ are faithful. If
$\psi\varphi=m\chi$ where $m$ is a positive integer and $\chi\in\Irr(G)$ then
$\psi$ and $\varphi$ are fully ramified with respect to $\Z(G)$.
\end{conA}

Other ways to state the conclusion of this conjecture are that $\varphi,\psi$
and $\chi$ vanish on $G-\Z(G)$ or that
$\varphi(1)=\psi(1)=\chi(1)=|G:\Z(G)|^{1/2}$ (by Problem 6.3 of \cite{isa}).
In particular, if $m=1$, these equalities yield $\varphi(1)=1$ and since it is
faithful, we deduce that $G$ is cyclic. So Conjecture A is indeed a strong form
of Isaacs conjecture.

Among other results, Isaacs proved that a counterexample to his conjecture has
Fitting height at least $4$ (see Theorem A of \cite{isa00}). We can prove
Conjecture A for nilpotent groups.

\begin{thmB}
Conjecture A holds for $p$-groups.
\end{thmB}

Using Theorem B we can prove Conjecture A for $p$-special characters (see
\cite{gaj} for their definition and basic properties).

\begin{thmC}
Let $G$ be a $p$-solvable group and suppose that the product of two faithful
$p$-special characters is a multiple of a $p$-special character. Then $G$ is a
$p$-group.
\end{thmC}

Theorem C is an easy consequence of the following elementary, but
perhaps surprising, result.

\begin{thmD}
Let $\varphi$ be a faithful  character of a finite group $G$ and
assume that $\psi\in\Irr(G)$. Write $\varphi\psi=m\Delta$, where $\Delta$ is a  character of $G$. If
$\Delta(1)\leq  \psi(1)$, then all three $\varphi, \psi$ and  $\Delta$  vanish on  $G-\Z(G)$. 
If in addition, $\varphi$ is irreducible, then $\varphi(1) = \psi(1) = \Delta(1)$.
\end{thmD}

Observe that Theorem D proves Conjecture A in the case that $\chi(1)$ is no larger than  at least one of $\varphi(1)$ and $\psi(1)$.

We thank E. C. Dade and I. M. Isaacs for many useful conversations. We also thank the referee for pointing out Lemma 3.1 to us, and thus making 
the proof of Theorem D much easier than our original one.
  The work of the third author was done while he was visiting the universities
of Crete and Wisconsin, Madison. He thanks both Departments
for their hospitality.

\section{Proof of Theorem B}

We begin work toward a proof of Theorem B. We need two elementary lemmas.

\begin{lem}\label{lem.1}
Let $\chi \in \Irr(G)$, where $G$ is a $p$-group.
Suppose $Z \subseteq Y \nor  G$, where $Z \nor G$ and $|Y:Z| = p$. If
$Z \sbs \Z( \chi) $ and $Y \nsubseteq  \Z(\chi)$, then $\chi$ vanishes
on $Y - Z$.
\end{lem}

\begin{proof}
Let $\lambda$ be the unique (linear) irreducible
constituent of $\chi_Z$. Then every irreducible constituent of $\chi_Y$ is an
extension of $\lambda$, and in particular is linear.  Because $Y \nsubseteq 
\Z(\chi)$, the number of distinct
linear constituents of $\chi_Y$ is a power of $p$ exceeding $1$, and
so is at least $p$. It follows that the irreducible constituents of $\chi_Y$
are all of the extensions of $\lambda$, and they all occur with equal
multiplicity, as $Y \unlhd G$. Since the sum of these extensions is
$\lambda^Y$, that sum
vanishes on $Y - Z$ and the result follows.
\end{proof}

\begin{lem}\label{lem2}
 Let $\epsilon$ and $\delta$ be $p$th roots of unity,
where $p$ is an odd prime. If $\delta \ne 1$, then 
$$
\left|\sum_{i=0}^{p-1} \epsilon^i \delta^{i(i-1)/2}\right| = \sqrt p \,.
$$
\end{lem}

\begin{proof}
Write $\epsilon = \delta^k$, where $0 \le k < p$. Then the
$i$th term of the sum is $\delta^{ki + i(i-1)/2}$. Since $p \ne 2$, we can
write this in the form $\delta^{ai^2 + bi}$ for suitable integer constants $a$
and $b$, where $1 \le a < p$ and $0 \le b < p$. Let $\tau = \delta^a$. The
$i$th term of our sum is then $\tau^{i^2 + 2ci}$ for some constant $c$. If
we multiply the sum by $\tau^{c^2}$, the $i$th term becomes
$\tau^{(i+c)^2}$. Since $i+c$ runs over the same set of values (mod $p$) as
$i$, we can rewrite our sum as $\sum \tau^{i^2}$. This is the well known
Gauss
sum, with absolute value $\sqrt p$ (see, for instance, p. 84 of
\cite{lan}).
\end{proof}

The next result is Theorem B.

\begin{teo}
Assume that $G$ is a  finite $p$-group,  for some prime $p$.
Assume further that $\varphi$ and  $\psi$ are faithful irreducible  characters
of $G$  
whose product is a multiple of an irreducible character.
Then $\varphi$ and $\psi$ vanish on $G - \Z(G)$.
\end{teo}

\begin{proof}
Let $\varphi \psi= m \chi$, for some positive integer $m$ and an irreducible
character 
$\chi$ of $G$.  
We argue by induction on $|G|$. So  assume that $G$ is a minimal
counterexample.   
Clearly $G$ is not abelian.  So the center  $\Z(G)$ of $G$  is a cyclic  proper
subgroup of $G$, since $G$ 
has a faithful irreducible character $\varphi$.

\setcounter{step}{0}

\begin{step} $G$ has an elementary abelian normal subgroup of order $p^2$.
\end{step}
Assume  that every normal abelian subgroup of $G$ is cyclic.
Then 4.3 of \cite{suz} yields that $G$ is dihedral or semidihedral of order
$\geq16$ or (generalized) quaternion. 
If $G\cong Q_8$, the result is clear. Thus, we may assume that $|G|\geq16$. 
Since $\varphi$ and $\psi$ lie over the unique non-principal irreducible
character 
of $\Z(G)$, $\chi$ lies over $1_{\Z(G)}$. Also, it is clear that
$\varphi(1)=\psi(1)=2$ and they vanish 
on $G-G'$. It follows that $\chi$ is not linear, i.e, $\chi(1)=2$ and $m=2$. 
Pick $x\in\Z_2(G)-\Z(G)$. We have that
$$
4=2|\chi(x)|=|\varphi(x)||\psi(x)|<4,
$$
because $x\in\Z(\chi)$ but $x\not\in\Z(\varphi)$. This contradiction
proves Step 1.

We fix an elementary abelian normal subgroup $A$ of $G$ of order $p^2$.  
Then $Z= A \cap \Z(G)$ is the cyclic group  $\Omega_1(\Z(G))$ of order
$p$. Furthermore  
$A = K  \times Z$, for some subgroup $K$ of $G$ of order $p$. 
Put $C=\C_G(A)=\C_G(K)$. 
Note that $|G:C|=p$.

\begin{step}
$\varphi_C$ and $\psi_C$ are reducible and each has a
unique irreducible constituent with kernel containing $K$.
\end{step}

The center $\Z(C)$ certainly contains $A$ and thus is not cyclic.
Hence $C$ has no irreducible faithful character.
Therefore   $\varphi_C$ and $\psi_C$ reduce. Because $C$ has index $p$ in $G$, 
both  $\varphi_C$ and $\psi_C$ equal  the sum of $p$ distinct  irreducible
constituents 
that form a single $G$-orbit.
Let $\varphi_1$ be an irreducible constituent of $\varphi_C$. 
Note that $A \cap \Ker{\varphi_1}$ is nontrivial since $A \sbs \Z( C)$ and $A$
is noncyclic.
 Also, this intersection does not contain $Z$ since $Z \nor G$ and $\varphi$ is
faithful. 
It follows that $A \cap \Ker{\varphi_1}$ is one of the $p$ subgroups of order
$p$ in $A$ other than $Z$, and we note that these subgroups form a
$G$-conjugacy class. We can thus replace $\varphi_1$ by a $G$-conjugate
and assume that $K = A \cap \Ker{\varphi_1}$. Similarly, $\psi_C$ is
reducible and has an irreducible constituent, say $\psi_1$, with kernel
containing $K$.

Since $K = A \cap \Ker{\varphi_1}$ and $K$ has $p$ distinct conjugates in
$G$, it follows that the subgroups $A \cap \Ker{\varphi_i}$ are distinct as
$\varphi_i$ runs over the $p$ irreducible constituents of $\varphi_C$. This
establishes the uniqueness for $\varphi_1$ and a similar argument works for
$\psi_1$.

We now  fix  irreducible constituents $\varphi_1$ and $\psi_1$ of $\varphi_C$
and $\psi_C$  respectively, such that $K \sbs \Ker{\varphi_1}$ and $K \sbs
\Ker{\psi_1}$.

\begin{step}
$\varphi$ and $\psi$ vanish on $G - C$ and $\chi$ is
faithful. Also, $\chi_C$ is reducible and $\varphi_1 \psi_1 = (m/p)\chi_1$,
where $\chi_1$ is the unique irreducible constituent of $\chi_C$ with
kernel containing $K$.
\end{step}

According to Clifford's theorem,   $\varphi = \varphi_1^G$ and $\psi =
\psi_1^G$, and thus 
 $\varphi, \psi$  vanish on  $G - C$.  Hence $\chi$ vanishes on $G - C$. It
follows
that $\chi_C$ is reducible, and thus is the sum of $p$ distinct irreducible
constituents. 
Since $K$ is in the kernel of both $\varphi_1$ and
$\psi_1$, we see that $\varphi_1\psi_1$ is a sum (with multiplicities) of
irreducible constituents of $\chi_C$ having $K$ in their kernel. 

Let $\psi_2$ be an irreducible constituent of $\psi_C$ different from
$\psi_1$. Then $K$ is not in the kernel of $\psi_2$, and so it is not in the
kernel of $\varphi_1\psi_2$. It follows that $K$ is in the kernel of some
irreducible constituent $\chi_1$ of $\chi_C$ but $K$ is not in the kernel of
all of the conjugates of $\chi_1$.

If $Z \sbs \Ker\chi$ then $A = ZK \sbs \Ker\chi_1$, and since $A \nor G$,
we see that $A$ is contained in the kernel of every irreducible
constituent of $\chi_C$, which is not the case. Thus $Z \nsubseteq  \Ker\chi$. 
This,  along with the fact that $\Z(G)$ is cyclic, implies that $\chi$ is
faithful.

Therefore, the same argument we gave in Step 2 for $\varphi$, implies 
  that $\chi_1$ is the unique irreducible constituent
of $\chi_C$ with kernel containing $K$. It follows that
$\varphi_1\psi_1 = m_1\chi_1$ for some integer $m_1$. Comparison of degrees
yields $(\varphi(1)/p)(\psi(1)/p) = m_1(\chi(1)/p)$. Since
$\varphi(1)\psi(1) = m\chi(1)$, we deduce that $m_1 = m/p$.

\begin{step} $p \ne 2$.

\end{step}

Otherwise $|Z| = 2$  and $Z$ has a unique nonprincipal irreducible character. 
In this case, both $\varphi$ and $\psi$ lie above this nonprincipal character. 
Hence $Z \sbs \Ker\varphi\psi$. Then $Z \sbs\Ker\chi$, which is not the
case.

Let $V/K = \Z(C/K)$ and write $Y = A\Z( G)$. Note that $Y \nor G$ and
that $Y \sbs \Z( C) \sbs V$.

\medskip

\begin{step} $V > Y$.

\end{step}

Note that $Y = K\Z( G) $ and assume that $V = K\Z( G) $. Let
$K_1 = \Ker\varphi_1$. If $K_1 > K$, then $(K_1/K) \cap \Z(C/K) > 1$,
and thus $K_1 \cap K\Z( G) = K_1 \cap V > K$. It follows that $K_1 \cap
\Z( G )> 1$, and thus $Z \sbs K_1$ as $\Z(G)$ is cyclic.
 This is not the case, however, since $Z\nsubseteq  \Ker\varphi_1$. We conclude
that $K_1 = K$.

Similarly we show that $\Ker\psi_1 = K$. Hence $\varphi_1$, $\psi_1$ are
inflated from unique 
faithful  characters $\bar{\varphi}_1$ and $\bar{\psi}_1$, respectively
of the
factor group 
$C/K$. In addition, $\chi_1$ is also inflated from a unique character
$\bar{\chi}_1$ of $C/K$ and satisfies
$\bar{\varphi}_1 \bar{\psi}_1 = m_1 \bar{\chi}_1$. By the minimality of
$G$, we
conclude that 
$\varphi_1$ and $\psi_1$  vanish on $C - V$. 

In this situation, where $V = Y$, we see that $V \nor G$, and thus all
irreducible constituents of $\varphi_C$ vanish on $C - V$. We conclude that
$\varphi$ vanishes on $G - V$. But $|V: \Z(G)| = p$ and $V \nsubseteq
\Z(\varphi)$. By Lemma \ref{lem.1}, therefore, $\varphi$ vanishes on $V - \Z(
G)$,
and hence on $G - \Z( G) $. Similarly, $\psi$ vanishes on $G - \Z(G)$,
and this is a contradiction since $G$ is a counterexample.
\medskip

\begin{step} $\Z(C )> Y$.

\end{step}

Certainly, $Y = A \Z( G) \sbs \Z(C)$ and we suppose that
equality occurs. Since $V > Y$, we can choose a subgroup $U$ such that
$Y \sbs U \sbs V$ and $|U:Y| = p$. We have $1 < [C,U] \sbs [C,V] \sbs K$,
and thus $[C,U] = K$. In particular, we see that $U \sbs \Z(\varphi_1)$
and so all values of $\varphi_1$ on $U$ are nonzero.

Now let $\varphi_2$ be any irreducible constituent of $\varphi_C$ other
than $\varphi_1$. We argue $U \not\sbs \Z(\varphi_2)$ since otherwise,
$K = [C,U] \sbs \Ker\varphi_2$, which is not the case. But $Y \sbs \Z( C)$,
and $|U:Y| = p$, so Lemma \ref{lem.1} implies that $\varphi_2$ vanishes
on $U - Y$. Since $\varphi_2$ is an arbitrary constituent of $\varphi_C$
other than $\varphi_1$, it follows that if $u \in U - Y$, then
$$
\varphi(u) =  \sum_{i=1}^p \varphi_i(u) = \varphi_1(u) \ne 0.
$$ 
Similarly, $\psi(u) = \psi_1(u) \ne 0$
and also $\chi(u) = \chi_1(u) \ne 0$.
We now have
$$
m\chi_1(u) = m\chi(u) = \varphi(u)\psi(u) = \varphi_1(u)\psi_1(u) =
(m/p)\chi_1(u)
$$
and this is a contradiction.
\medskip

We choose $W \nor G$ with $Y \sbs W \sbs \Z( C) $ and $|W:Y| = p$. We also fix
elements $g \in G - C$ and $w \in W - Y$ and we write $[w,g] = a$ and $[a,g] =
z$. Then we can show
\medskip

\begin{step} $1 \ne a \in A $, $1 \ne z \in Z$ and $w^{g^i} =wa^iz^{i(i-1)/2}$.
\end{step}

Since $|W:Y| = p$, we have $W/Y \sbs \Z(G/Y)$, and thus $[W, G] \leq Y$. 
 Hence 
 $a = [w,g]$ is an element of $Y$ and thus $a^p \in \Z(G)$. 
Also, $|Y:\Z( G) | = p$, and  similarly we get $z \in \Z( G )$. Then
$a^p = (a^p )^g = (az)^p = a^p z^p$.  We conclude that   
$z \in Z= \Omega_1(\Z(G))$.

Because $w^g = w a$ and $a^g = a z$ we can easily calculate that 
$w^{g^i} = wa^iz^{i(i-1)/2}$ for integers $i$ with $1 \leq i  \leq p$. Note
that  
 $g^p \in C$ while $w  \in W \leq \Z(C)$. So  $w^{g^p} = w$. 
Also, since $p \ne 2$ and $z^p = 1$, we see that $z^{p(p-1)/2} = 1$. It 
follows
that $w = w^{ g^p}= w a^p$ and so $a^p = 1$. 
Since $a \in Y$ and $A = \Omega_1(Y)$, we have $a \in A$,
as wanted.

Finally, we must show that $ a \ne 1$ and  $z \ne 1$. 
If  $a =1 $ then $w \in \Z(C)$ is centralized by $g$, and thus 
$w \in \Z(G)$ contradicting the way $w$ was picked.
Also if $z =1$,  then    $g$ centralizes
$a \in \Z(C)$,   and thus  $a \in \Z( G) $. Hence $a \in A \cap \Z(G) = Z$.
  Note that since $A$ is not central in $G$  we have  $1 < [A,G] \nor G$. 
It follows that $[A,G] = Z$ and thus $[Y,G] = Z$. 
Hence $[W,g] \sbs Z$ since $W = Y \langle w \rangle$.
 But $W$ is abelian, and it
follows that $|W: \C_W( g) | \le |Z| = p$. This is a contradiction,
however
since $\C_W(g) = \Z( G) $ has index $p^2$ in $W$.

\medskip

\begin{step}

We have a contradiction.

\end{step}

Since $W \sbs \Z( C) $, there exists a linear character $\alpha \in \Lin(W)$
such that 
$(\varphi_1)_W = \varphi_1(1) \alpha$. Furthermore,  as $A \subseteq W$
we can write 
$\alpha(a) = \epsilon $ and $\alpha(z) = \delta$,
where $\epsilon$ and $\delta$ are $p$th roots of unity and $\delta \ne 1$
since $z \ne 1$ and $\varphi$ is faithful. 
We see now that
$$
\varphi(w) = \sum_{i =0}^{p-1} \varphi_1({w^{g^i}} )=\sum_{i=0}^{p-1}
\varphi_1(1) \alpha(w) \alpha(a^i) \alpha( z^{i(i-1)/2}) = 
 \varphi_1(w) A \,,
$$
where $A = \sum_{i=0}^{p-1} \epsilon^i\delta^{i(i-1)/2}$. By Lemma \ref{lem2},
therefore,
we have $|A| = \sqrt p$. We also have similar formulas
$\psi(w) = \psi_1(w) B$ and $\chi(w) = \chi_1(w) D$, where
$|B| = |D| = \sqrt p$. Also $\chi_1(w) \ne 0$ since $w \in \Z( C )$.

We have
$$
m \chi_1(w) D = m \chi(w) = \varphi(w)\psi(w) =
\varphi_1(w)\psi_1(w) AB = (m/p)\chi_1(w) AB \,,
$$
and thus $AB = pD$. But this is not consistent with
$|A| = |B| = |D| = \sqrt p$, and the proof is complete.
\end{proof}

\section{Proof of Theorems C and D}

In order to prove Theorem D,
we need the following  lemma.

\begin{lem}
\label{cs}
Let $\varphi \psi = m \Delta$, where $\Delta(1) \leq \psi(1)$, where $\psi$ is an irreducible character of $G$. 
Then \\
{\rm (a)}  $\Delta$ is irreducible. \\
{\rm (b)} $\Delta(1) = \psi(1)$.\\
{ \rm (c)} $\varphi$ and $\psi$ both vanish on $G - \Z(\varphi)$.
\end{lem}

\begin{proof}
We have $m \geq \varphi(1)$. If $\varphi(x) = 0 $ then $\Delta(x) = 0$. Otherwise, 
\begin{equation*}
|\psi(x) | = \frac{m |\Delta(x)|} {|\varphi(x)|} \geq \frac{\varphi(1)}{|\varphi(x)|} |\Delta(x)| \geq |\Delta(x)|.
\end{equation*}
Therefore in all cases we get $|\psi(x)| \geq |\Delta(x)|$.  So $1 = [ \psi, \psi]  \geq [\Delta, \Delta]$. Because $\Delta$ is a character, 
we have equality everywhere. Therefore $|\Delta(x) | = |\psi(x)|$ for all $x \in G$.  Hence (a) and (b) follow.

If $\varphi(x) \ne 0$ then $|\varphi(x)| = \varphi(1)$, and therefore $x \in \Z(\varphi)$. Hence for every  $x \in G - \Z(\varphi)$ we have $\varphi(x) = 0$, and thus 
$\Delta(x) = 0$.  Because $|\psi(x) |= |\Delta(x)|$ for all $x$, we also get $\psi(x) = 0$, for $x \in G - \Z(\varphi)$. This proves (c).
\end{proof}

\begin{proof}[Proof of Theorem D]
If $\varphi$ is faithful , then $\Z(\varphi)  \subseteq  \Z(G)$. According to the above lemma, the irreducible character vanishes outside the $\Z(\varphi)$. 
Therefore, $\Z(\varphi) = \Z(G)$. In addition, the above lemma implies that all three $\varphi$ , $\psi$ and $\Delta$ vanish on $G - \Z(G)$.
Hence Theorem D follows. 
\end{proof}

Theorem C is an immediate consequence of the following result.

\begin{teo}\label{tt}

Let $G$ be a finite group and suppose that $\varphi,\psi\in\Irr(G)$ are
faithful and $\varphi\psi=m\chi$ with 
$\chi\in\Irr(G)$. Let $P$ be a $p$-Sylow  subgroup of $G$ for some prime $p$. If 
$\varphi_P, \psi_P$ and  $\chi_P$ are irreducible then $G$ is nilpotent.
\end{teo}

\begin{proof}

Assume that all three restrictions to $P$ are
irreducible. We can apply Theorem B and deduce 
that $\chi(1)=\varphi(1)=\psi(1)$.  Furthermore, because they restrict irreducibly to $P$, 
their degrees are $p$-numbers. 
By Theorem D, (since $\chi(1) \leq  \psi(1)$), we have that $\chi$ vanishes on $G -  \Z(G)$.
On the other hand, since the degree of the  irreducible character $\chi$ is not divisible by any prime $q \ne p$, 
$\chi$ never vanishes on any element of order a power of $q$ (see Lemma 2.1 in \cite{na}). Hence $Q \subseteq \Z(G)$, 
for any  $q$-Sylow subgroup $Q$ of $G$, with $q \ne p$.
This is enough to prove the theorem.
\end{proof}

Finally, we complete the proof of Theorem C.

\begin{proof}[Proof of Theorem C] 
Note that the restriction of $p$-special characters to a Sylow $p$-subgroup is
irreducible (see \cite{gaj}).
 Hence  Theorem \ref{tt} implies that $G$ is nilpotent.  But $G$ has  a faithful $p$-special character, therefore 
it has to be a $p$-group.
\end{proof}


\begin{thebibliography}{99}



\bibitem{gaj} D. Gajendragadkar, A characteristic class of 
characters of finite $\pi$-separable groups, {\it J. Algebra} {\bf 59} 
(1979), 237--259.



\bibitem{isa}
M. Isaacs, ``Character Theory of Finite Groups",
Dover, New York, 1994.



\bibitem{isa00} 
M. Isaacs, Irreducible products of characters, {\it J. Algebra} {\bf 223}
(2000), 630--646.


\bibitem{lan}
S. Lang, ``Algebraic Number Theory", Springer-Verlag, New York, 1986.

\bibitem{na}
G. Navarro, 
Zeros of primitive  characters, {\it J. Algebra} {\bf 221} (1999),
644--650. 


\bibitem{suz}
M. Suzuki, ``Group Theory II", Springer-Verlag, New York, 1986.



\end{thebibliography}
\end{document}